\documentclass[12pt]{amsart}

\usepackage{amsfonts}
\usepackage{amssymb}
\usepackage{eucal}

\begin{document} 

\textheight     230.0mm
\textwidth      131.0mm
\oddsidemargin    0.0mm
\evensidemargin   0.0mm
\topmargin        0.0mm
\arraycolsep      2.0pt

\font\smss=cmss10 scaled\magstep1
\font\eightsans=cmss8
\font\smallit=ptmri at 10pt
\font\sii=ptmri at 8pt
\font\eightrm=cmr8
\font\eightbf=cmbx8
\font\titlebf=cmbx10 scaled\magstep2
\font\eb=cmbx8
\font\eightit=cmti8
\font\eighti=cmmi8
\font\eightsy=cmsy8
\font\eightsmc=cmcsc8
\font\eighttt=cmtt8
\font\sixrm=cmr6
\font\sixi=cmmi6
\font\sixsy=cmsy6
\font\tenex=cmex10
   \def\eightpoint{%
   \def\rm{\fam0\eightrm}\def\bf{\fam\bffam\eightbf}%
  \def\it{\fam\itfam\eightit}\def\smc{\eightsmc}\def\tt{\eighttt}\baselineskip=10pt\rm%
        \textfont0=\eightrm \scriptfont0=\sixrm
        \textfont1=\eighti \scriptfont1=\sixi
        \textfont2=\eightsy \scriptfont2=\sixsy
\textfont3=\tenex \scriptfont3=\tenex
}
\def\hs{\hskip.7pt}
\def\hh{\hskip.4pt}
\def\nh{\hskip-.7pt}
\def\nnh{\hskip-1pt}
\def\bbR{\text{\bf R}}
\def\er{r}
\def\bs{\varSigma}
\def\bx{L}
\def\cy{{y}}
\def\cj{{c}}
\def\tb{T\hskip-.3pt\bs}
\def\tab{{T\hskip.2pt^*\hskip-2.3pt\bs}}
\def\tacb{{T_\cy^*\hskip-1.9pt\bs}}
\def\txm{{T_x\hskip-.3ptM}}
\def\txmy{{T_x\hskip-.3ptM_\cy}}
\def\vy{V_{\hskip-1.5pt\cy}}
\def\tm{{T\hskip-.3ptM}}
\def\vt{\mathcal{P}}
\def\w{\vt^\perp}
\def\vd{\vt\hh'}
\def\vdx{\vd{}\hskip-4.5pt_x}
\def\cc{\mathcal{C}}
\def\dd{\mathcal{D}}
\def\ee{\mathcal{E}}
\def\xz{\mathcal{X}}
\def\pz{z}
\def\gm{h}
\def\hj{\gamma}
\def\lgl{\langle}
\def\rgl{\rangle}
\def\lr{\langle\,,\rangle}
\def\vs{vector space}
\def\rvs{real vector space}
\def\vf{vector field}
\def\dn{distribution}
\def\inv{-in\-var\-i\-ant}
\def\diml{-di\-men\-sion\-al}
\def\prl{-par\-al\-lel}
\def\dbly{-dif\-fer\-en\-ti\-a\-bly}
\def\kx{\alpha}
\def\mf{manifold}
\def\bmf{base manifold}
\def\pb{-prin\-ci\-pal bundle}
\def\pbd{N}
\def\bd{bundle}
\def\bp{bundle projection}
\def\prc{pseu\-do\hs-Riem\-ann\-i\-an metric}
\def\prd{pseu\-\hbox{do\hs-}Riem\-ann\-i\-an manifold}
\def\npd{null parallel distribution}
\def\pj{-pro\-ject\-a\-ble}
\def\lcc{Le\-vi-Ci\-vi\-ta connection}
\def\vb{vector bundle}
\def\vbm{vec\-tor-bun\-dle morphism}
\def\afb{affine bundle}
\def\aff{A}
\def\kerd{\text{\rm Ker}\hskip2.7ptd}
\def\ts{total space}
\def\pmb{\pi}
\def\pqb{\hs\text{\rm p}}

\newtheorem{thm}{Theorem}[section] 
\newtheorem{prop}[thm]{Proposition} 
\newtheorem{lem}[thm]{Lemma} 
\newtheorem{cor}[thm]{Corollary} 
  
\theoremstyle{definition} 
  
\newtheorem{defn}[thm]{Definition} 
\newtheorem{notation}[thm]{Notation} 
\newtheorem{example}[thm]{Example} 
\newtheorem{conj}[thm]{Conjecture} 
\newtheorem{prob}[thm]{Problem} 
  
\theoremstyle{remark} 
  
\newtheorem{rem}[thm]{Remark}

\baselineskip=15.5pt
\parskip=3.5pt
\renewcommand{\thesection}{\Roman{section}}
\renewcommand{\thethm}{\arabic{section}.\arabic{thm}}

\ 

\ 

\

\ 

\ 

\noindent{\titlebf Walker\hskip.6pt's theorem without coordinates}

\ 

\vskip8pt

\hskip22pt\parbox[l]{4.5246in}{
\noindent{\smss Andrzej Derdzinski}

\noindent{\smallit Department of Mathematics, The Ohio State University,\\ 
Columbus, OH 43210, USA}

\noindent{\eightrm Electronic mail: andrzej@math.ohio-state.edu}

\ 

\noindent{\smss Witold Roter}

\noindent{\smallit Institute of Mathematics and Computer Science, Wroc\l aw 
University of Technology,\\
Wy\-brze\-\vbox{\hbox{\hskip1.4pt.\hskip-1.4pt}\vskip-8.2pt\hbox{\smallit z}}e 
Wys\-pia\'n\-skiego 27, 50-370 Wroc\l aw, Poland}

\noindent{\eightrm Electronic mail: roter@im.pwr.wroc.pl}
%
%

\ 

\noindent We provide a co\-or\-di\-nate-free version of the local 
classification, due to A.\hskip1.6ptG.\hskip2.3ptWalker [Quart.\ J.\ Math.\ 
Oxford (2) {\bf1}, 69 (1950)], of null parallel distributions on 
pseu\-\hbox{do\hs-}Riem\-ann\-i\-an manifolds. The underlying manifold is 
realized, locally, as the total space of a fibre bundle, each fibre of which 
is an affine principal bundle over a pseu\-\hbox{do\hs-}Riem\-ann\-i\-an 
manifold. All structures just named are naturally determined by the 
distribution and the metric, in contrast with the non\-ca\-non\-i\-cal choice 
of coordinates in the usual formulation of Walker's theorem.
}



\voffset=-28pt\hoffset=52pt  

\setcounter{section}{0}
\section{Introduction}\label{intr}
In 1950 A.\hskip1.6ptG.\hskip2.3ptWalker$^1$ described the local structure of 
all pseu\-do\hs- Riem\-ann\-i\-an manifolds with \npd s. The present paper 
provides a co\-or\-di\-nate-free version of Walker's theorem.

Many \hskip-.4ptauthors, \hskip-.4ptbeginning with \hskip-.9ptWalker 
himself,$^2$ have invoked \hskip-.9ptWalker's 1950 result, often to generalize 
it or derive other theorems from it. In our bibliography, which is by no means 
complete, \hbox{Refs.\ 3 -- 16} all belong to this category. They invariably 
cite Walker's result in its original, lo\-cal-co\-or\-di\-nate form 
(reproduced in the Appendix).

Such an approach, perfectly suited for the applications just mentioned, tends 
nevertheless to obscure the geometric meaning of Walker's theorem. In fact, 
Walker coordinates are far from unique; choosing them results in making 
non\-ca\-non\-i\-cal objects a part of the structure.

To keep the picture canonical, some authors$^{\hh3\hh,\hh5}$ replace a single 
Walker coordinate system by a whole maximal atlas of them. What we propose 
here, instead, is to use only ingredients such as fibre bundles, widely seen 
as more directly ``geometric'' than a coordinate atlas (even though one may 
ultimately need atlases to define them).

In our description, the co\-or\-di\-nate-in\-de\-pen\-dent content of Walker's 
theorem amounts to realizing the underlying manifold, locally, as a fibre 
bundle whose fibres are also bundles, namely, affine principal bundles over 
\prd s. The bundle structures are all naturally associated with the original 
\npd; the \dn\ and the metric can in turn be reconstructed from them.

\section{Preliminaries}\label{prel}
Throughout this paper, all manifolds, \bd s, sections, sub\bd s, connections, 
and mappings, including \bd\ morphisms, are assumed to be of class 
$\,C^\infty\nnh$. A \bd\ morphism may operate only between two \bd s with the 
same \bmf, and acts by identity on the base.

`\nnh A \bd' always means `a $\,C^\infty$ locally trivial \bd' and the same 
symbol, such as $\,M$, is used both for a given \bd\ and for its \ts; the \bd\ 
projection $\,M\to\bs\,$ onto the \bmf\ $\,\bs\,$ is denoted by $\,\pmb\,$ 
(or, sometimes, $\,\pqb$). We let $\,M_\cy$ stand for the fibre 
$\,\pmb^{-1}(\cy)$ over any $\,\cy\in\bs$, while $\,\hs\kerd\pmb\,$ is the 
vertical distribution treated as a \vb\ (namely, a sub\bd\ of the tangent \bd\ 
$\,\tm$).

For real \vb s $\,\mathcal{X},\mathcal{Y}\,$ over a manifold $\,\bs\,$ and a 
\rvs\ $\,V$ with $\,\dim V<\infty$, we denote by 
$\,\hs\text{\rm Hom}\hs(\mathcal{X},\mathcal{Y})\,$ the \vb\ over 
$\,\bs\,$ whose sections are \vbm s $\,\mathcal{X}\to\mathcal{Y}$, and by 
$\,\bs\times V$ the product \bd\ with the fibre $\,V\nnh$, the sections of 
which are functions $\,\bs\to V\nnh$. Thus, 
$\,\mathcal{X}\hh^*\nh=\,\text{\rm Hom}\hs(\mathcal{X},\bs\times\bbR)\,$ is 
the dual of $\,\mathcal{X}$.

We will say that a given fibrewise structure in a \bd\ $\,M\,$ over a \mf\ 
$\,\bs\,$ {\it depends\/ $\,C^\infty\nh$\dbly\ on\/} $\,\cy\in\bs$, or {\it 
varies\/ $\,C^\infty\nh$\dbly\ with\/} $\,\cy$, if suitable $\,C^\infty$ 
local trivializations of $\,M\,$ make the structure appear as constant (the 
same in each fibre).

The symbol $\,\nabla\,$ will be used for various connections in \vb s. Our 
sign convention about the curvature tensor $\,R=R^\nabla$ of a connection 
$\,\nabla\,$ in a \vb\ $\,\xz\,$ over a \mf\ $\,\bs\,$ is
\begin{equation}
R(u,v)\psi\hskip7pt=\hskip7pt\nabla_{\!v}\nabla_{\!u}\psi\,
-\,\nabla_{\!u}\nabla_{\!v}\psi\,+\,\nabla_{[u,v]}\psi\hs,
\label{cur}
\end{equation}
for sections $\,\psi\,$ of $\,\xz\,$ and \vf s $\,u,v\,$ tangent to 
$\,\bs$. By the Leibniz rule, when $\,\nabla$ is the \lcc\ of a \prc\ $\,g\,$ 
and $\,u,v,w\,$ are tangent \vf s, $\,2\hs\lgl\nabla_{\!w}v,u\rgl\,$ 
equals$^{17}$
\begin{equation}
d_w\lgl v,u\rgl+d_v\lgl w,u\rgl-d_u\lgl w,v\rgl+\lgl v,[u,w]\rgl
+\lgl u,[w,v]\rgl-\lgl w,[v,u]\rgl\hh,
\label{lcc}
\end{equation}
where $\,d_v$ is the directional derivative and $\,\lr\,$ stands for 
$\,g(\hskip3pt,\hskip1.8pt)$.
\begin{rem}\label{liebr}Let $\,\pmb:M\to\bs\,$ be a \bp. A \vf\ $\,w$ on the 
\ts\ $\,M\,$ is $\,\pmb$\pj\ onto the \bmf\ $\,\bs\,$ if and only if, for 
every vertical \vf\ $\,u\,$ on $\,M$, the Lie bracket $\,[w,u]\,$ is also 
vertical. This well-known fact is easily verified in local coordinates for 
$\,M\,$ which make $\,\pmb\,$ appear as a standard Euclidean projection.
\end{rem}

\section{Affine principal bundles}\label{afpb}
All principal \bd s discussed below have Abelian structure groups $\,G$, so 
one need not decide whether $\,G\,$ acts from the left or right.

Let $\,\pbd\hs$ be a $\,G$\pb\ over a \bmf\ $\,\bx$, where $\,G\,$ is an 
Abelian Lie group. By the $\,\pbd${\it-pro\-lon\-ga\-tion\/} of the tangent 
\bd\ $\,T\bx\,$ \hskip-.3ptwe mean the \vb\ $\,\mathcal{F}\,$ \hskip-.3ptover 
$\,\bx\,$ \hskip-.3ptwhose fibre $\,\mathcal{F}_{\hskip-1pt\cj}$ over 
$\,\cj\in\bx\,$ is the space of all $\,G$\inv\ \vf s tangent to $\,\pbd\hs$ 
along $\,\pbd_\cj$ (and defined just on $\,\pbd_\cj$), with $\,\pbd_\cj$ 
denoting, as usual, the fibre of $\,\pbd$ over $\,\cj$. A vector sub\bd\ 
$\,\mathcal{G}\subset\mathcal{F}\,$ now can be defined by requiring 
$\,\mathcal{G}_\cj$, for any $\,\cj\in\bx$, to consist of all $\,G$\inv\ \vf s 
defined just on $\,\pbd_\cj$ which are vertical (i.e., tangent to 
$\,\pbd_\cj$). Since each $\,\mathcal{G}_\cj$ is canonically isomorphic to the 
Lie algebra $\,\mathfrak{g}\,$ of $\,G$, the \vb\ $\,\mathcal{G}\,$ is 
naturally trivialized, that is, identified with the product \bd\ 
$\,\bx\times\hh\mathfrak{g}$. Therefore
\begin{equation}
\bx\times\hh\mathfrak{g}\,\,=\,\,\mathcal{G}\,\,\subset\,\,\mathcal{F}\hs.
\label{bge}
\end{equation}
The quotient \bd\ $\,\mathcal{F}/\mathcal{G}\,$ is in turn naturally 
isomorphic to $\,T\bx$, via the differential of the \bp\ $\,\pbd\to\bx$.

An {\it affine space\/} is a set $\,\aff\,$ with a simply transitive action on 
$\,\aff\,$ of the additive group of a vector space $\,V\nnh$. One calls $\,V$ 
the {\it vector space of translations\/} of the affine space $\,\aff$.

An {\it affine bundle\/} $\,M\,$ over a \mf\ $\,\bs\,$ is a \bd\ with fibres 
$\,M_\cy$, $\,\cy\in\bs$, carrying the structures of affine spaces whose \vs s 
$\,\mathcal{X}_\cy$ of translations form a \vb\ $\,\mathcal{X}\,$ over 
$\,\bs$, called the {\it associated vector bundle\/} of $\,M$. We also require 
the affine-space structure of $\,M_\cy$ to vary $\,C^\infty\nh$\dbly\ with 
$\,\cy\in\bs$, in the sense of \hbox{Sec.\ \ref{prel}}.

If, in addition, $\,\hs\mathcal{X}=\bs\times V\nnh$, that is, the associated 
\vb\ of $\,M\,$ happens to be a product \bd, then $\,M\,$ is also a $\,V$\pb, 
with the obvious action of the additive group of the vector space $\,V\nnh$. 
Such {\it affine principal bundles\/} are distinguished from arbitrary \afb s 
by having a structure group that, instead of general affine transformations of 
a model fibre, contains only translations.

\section{Partial metrics and extensions}\label{pmae}
Let $\,\cc$, $\,\dd\,$ and $\,\ee\,$ be real vector 
\bd s over a \mf\ $\,Q$. By an $\,\ee$-val\-ued {\it pairing\/} 
of $\,\hs\cc\,$ and $\,\dd\,$ we mean any \vbm\ 
$\,\beta:\cc\otimes\dd\to\ee$. This amounts to a 
$\,C^\infty$ assignment of a bilinear mapping 
$\,\beta(\pz):\cc_\pz\times\dd_\pz\to\ee_\pz$ to every 
$\,\pz\in Q$. An $\,\ee$-val\-ued {\it partial pairing\/} of 
$\,\hs\cc\,$ and $\,\dd\,$ consists, by definition, of two vector 
sub\bd s $\,\cc\hs'\subset\cc\,$ and 
$\,\dd\hh'\subset\dd$, of some codimensions $\,k\,$ and $\,l$, 
along with pairings 
$\,\hj:\cc\otimes\dd\hh'\to\ee\,$ and 
$\,\hj:\cc\hs'\otimes\dd\to\ee\,$ which coincide on 
the sub\bd\ $\,\cc\hs'\otimes\dd\hh'$ (and so may be 
represented by the same symbol $\,\hj$ without risk of ambiguity). 
One can obviously restrict a given pairing 
$\,\beta:\cc\otimes\dd\to\ee\,$ to 
$\,\cc\otimes\dd\hh'$ and 
$\,\cc\hs'\otimes\dd$, so that a partial pairing $\,\hj\,$ 
is obtained; we will then say that $\,\beta\,$ is a {\it to\-tal-pair\-ing 
extension\/} of $\,\hj$.
\begin{lem}\label{parex}For any fixed partial pairing $\,\hj$, with 
$\,\cc,\dd,\ee,\cc\hs'\nnh,\dd\hh'\nnh,k,l$ and\/ $\,Q\,$ as 
above, and with\/ $\,m\hs$ denoting the fibre dimension of $\,\ee$, the 
to\-tal-pair\-ing extensions of\/ $\,\hj\,$ coincide with sections of a 
specific \afb\ of fibre dimension\/ $\,klm\hs$ over $\,Q$, whose associated 
\vb\ is $\,\hs\text{\rm Hom}\hs(\cc/\cc\hs'\nh
\otimes\dd/\dd\hs'\nnh,\ee)$.
\end{lem}
\begin{proof}Our $\,\hj\,$ is nothing else than a \vbm\ 
$\,\xz\to\ee$, where  
$\,\xz\subset\cc\otimes\dd\,$ is the sub\bd\ spanned 
by $\,\cc\otimes\dd\hh'$ and 
$\,\cc\hs'\otimes\dd$. The \afb\ in question is the 
pre\-im\-age of the section $\,\hj\,$ under the (surjective) restriction 
morphism $\,\hs\text{\rm Hom}\hs(\cc\otimes\nh\dd,\ee)
\to\,\text{\rm Hom}\hs(\xz,\ee)$.
\end{proof}
As usual,$^5$ by a {\it pseu\-\hbox{do\hs-}Riem\-ann\-i\-an fibre metric\/} 
$\,g\,$ in a \vb\ $\,\mathcal{T}\,$ over a \mf\ $\,M\,$ we mean any family of 
nondegenerate symmetric bilinear forms $\,g(x)\,$ in the fibres 
$\,\mathcal{T}_x$ that constitutes a $\,C^\infty$ section of the symmetric 
power $\,(\mathcal{T}^*){}^{\odot2}\nnh$. Equivalently, such $\,g\,$ is a 
pairing of $\,\mathcal{T}$ and $\,\mathcal{T}\,$ valued in the product \bd\ 
$\,M\times\bbR\hs$, symmetric and non\-de\-gen\-er\-ate at every point of 
$\,M$.

Let $\,\mathcal{T}\,$ again be a \vb\ over a \mf\ $\,M$. We define a {\it 
partial fibre metric\/} in $\,\mathcal{T}\,$ to be a triple 
$\,(\vt,\vd\nnh,\kx)\,$ formed by vector sub\bd s $\,\vt\,$ and $\,\vd$ 
of $\,\hs\mathcal{T}\,$ along with a pairing 
$\,\kx:\vd\otimes\mathcal{T}\to M\times\bbR\hs$, valued in the product 
\bd\ $\,M\times\bbR\hs$, such that
\begin{enumerate}
  \def\theenumi{{\rm\roman{enumi}}}
\item $\mathcal{T},\hs\vt\,$ and $\,\vd$ are of fibre dimensions 
$\,n,\hs\er\,$ and, respectively, $\,n-\er$ for some $\,n,\hs\er\,$ with 
$\,0\le\er\le n/2$, while $\,\vt\subset\vd\nnh$,
\item at every $\,x\in M\,$ the bilinear mapping 
$\,\kx(x):\vdx\times\mathcal{T}_x\to\bbR\,$ has the rank $\,n-\er$, its 
restriction to $\,\vdx\times\vdx$ is symmetric, and its restriction to 
$\,\vdx\times\vt_x$ equals $\,0$.
\end{enumerate}
By a {\it to\-tal-met\-ric extension\/} of $\,(\vt,\vd\nnh,\kx)\,$ we then 
mean any \hbox{pseu\-do\hs-} Riem\-ann\-i\-an fibre metric in 
$\,\mathcal{T}\,$ whose restriction to $\,\vd\nh\otimes\mathcal{T}\,$ is 
$\,\kx$.
\begin{lem}\label{parfm}The to\-tal-met\-ric extensions\/ $\,g\hs$ of any 
partial fibre metric $\,(\vt,\vd\nnh,\kx)$, with $\,\er,\hs M\,$ as above, 
coincide with the sections of a specific \afb\ of fibre dimension $\,r(r+1)/2$ 
over $\,M$. For every such $\,g$ the sub\bd\ $\,\vt\hs$ is $\,g$-null and\/ 
$\,\vd\nh$ is its $\,g$-or\-thog\-o\-nal complement.
\end{lem}
\begin{proof}For any fixed point $\,x\in M\nh$, let us choose a basis 
$\,e_1,\dots,e_n$ of $\,\mathcal{T}_x$ such that $\,e_1,\dots,e_r\in\vt_x$ and 
$\,e_1,\dots,e_{n-r}\in\vdx$. The matrix of $\,g(x)$, for any to\-tal-met\-ric 
extension $\,g\,$ of our partial fibre metric, then is the matrix appearing in 
Walker's original theorem (see the Appendix), with $\,\det A\ne0$, and with 
the two occurrences of $\,I\,$ replaced by some nonsingular $\,r\times r\hs$ 
matrix $\,C\hs$ and its transpose $\,C\hh'\nnh$. The sub-matrices $\,A,H,C$ 
(and $\,H'\nnh,C\hh'$) are prescribed, while the freedom in choosing 
$\,g(x)\,$ is represented by an arbitrary symmetric $\,r\times r\hs$ matrix 
$\,B$.
\end{proof}

\section{Walker's theorem}\label{walk}
Suppose that the following data are given.
\begin{enumerate}
  \def\theenumi{{\rm\alph{enumi}}}
\item Integers $\,n\,$ and $\,\er\,$ with $\,0\le\er\le n/2$.
\item An $\,r$\diml\ \mf\ $\,\bs$.
\item A \bd\ over $\,\bs\,$ with some \ts\ $\,M\hskip-.3pt$, whose 
every fibre $\,M_\cy$, $\,\cy\in\bs$, is a $\,\tacb$-prin\-ci\-pal 
\bd\ over a $\,(n-2\er)$\diml\ \mf\ $\,Q_\cy$. (Cf.\ the last paragraph of 
\hbox{Sec.\ \ref{afpb}}.)
\item A \prc\ $\,\gm_\cy$ on each $\,Q_\cy$, $\,\cy\in\bs$.
\end{enumerate}
We assume that all $\,\cy$-dependent objects in (c) -- (d), including the 
prin\-ci\-pal-\bd\ structure, vary $\,C^\infty\nh$\dbly\ with $\,\cy\in\bs\,$ 
(in the sense of \hbox{Sec.\ \ref{prel}}) and, in particular, the $\,Q_\cy$ 
are the fibres of a \bd\ over $\,\bs$ with a \ts\ $\,Q\,$ of dimension 
$\,n-\er$. When $\,\er=n/2$, each $\,\gm_\cy$ is the ``zero metric'' on the 
discrete space $\,Q_\cy$, cf.\ \hbox{Sec.\ \ref{midi}}.

Let $\,\mathcal{F}\,$ be the \vb\ over $\,Q\,$ whose restriction to $\,Q_\cy$, 
for each $\,\cy\in\bs$, is the $\,M_\cy$-pro\-lon\-ga\-tion of the tangent 
\bd\ $\,TQ_\cy$ (see \hbox{Sec.\ \ref{afpb}}) for the $\,\tacb$\pb\ $\,M_\cy$ 
over $\,Q_\cy$. Relation (\ref{bge}) now yields 
$\,\pqb^*(\tab)\subset\mathcal{F}$, where $\,\pqb:Q\to\bs\,$ denotes the \bp. 
In other words, $\,\pqb^*(\tab)\,$ may be treated as a vector sub\bd\ of 
$\,\mathcal{F}$.

Furthermore, the quotient-\bd\ identification following formula (\ref{bge}) 
yields $\,\mathcal{F}/\pqb^*(\tab)=\,\kerd\pqb\,$ (the vertical sub\bd\ of 
$\,TQ$, for the projection $\,\pqb:Q\to\bs$). 

We define a partial pairing $\,\hj\,$ of $\,\mathcal{F}\,$ and $\,TQ\,$ valued 
in the product \bd\ $\,Q\times\nh\bbR\hs$, as in \hbox{Sec.\ \ref{pmae}}, for 
our $\,Q\,$ along with $\,\cc=\mathcal{F}$, $\,\dd=TQ$, 
$\,\ee=Q\times\nh\bbR\hs$, $\,\cc\hs'=\pqb^*(\tab)\,$ and 
$\,\dd\hh'=\,\kerd\pqb$. Namely, given $\,\pz\in Q$, 
we set $\,\hj(\xi,\zeta)=\xi(d\pqb_\pz\zeta)\,$ for 
$\,\xi\in\tacb=[\pqb^*(\tab)]_\pz$ and $\,\zeta\in T_\pz Q$, with 
$\,\cy=\pqb(\pz)\in\bs$, as well as 
$\,\hj(u,\psi)=\gm_\cy([u],\psi)\,$ for 
$\,u\in\mathcal{F}_{\hskip-1pt\pz}$ and $\,\psi\in\,\kerd\pqb_\pz$, where 
$\,u\mapsto[u]\,$ denotes the surjective \vbm\ 
$\,\mathcal{F}\to\hs\kerd\pqb\,$ with the kernel $\,\pqb^*(\tab)$.

Our construction has two steps involving arbitrary choices:
\vskip3pt
\noindent{\bf Step 1. }{\it We choose 
$\,\beta:\mathcal{F}\otimes TQ\to Q\times\bbR\,$ to be any to\-tal-pair\-ing 
extension of $\,\hj$.}
\vskip3pt
\noindent According to Lemma~\ref{parex}, such $\,\beta\,$ is just an 
arbitrary section of an \afb\ of fibre dimension $\,(n-2r)r\,$ over $\,Q$. For 
the meaning of the above discussion in Walker's original language, see the 
Appendix.

The remainder of our construction proceeds as follows. Using $\,\beta$, we 
define a partial metric $\,(\vt,\vd\nnh,\kx)\,$ in the tangent \bd\ $\,\tm$. 
Specifically, $\,\mathcal{T},\vt,\vd$ and $\,n,r\,$ with the properties listed 
in (i) -- (ii) of \hbox{Sec.\ \ref{pmae}} are chosen so that 
$\,\mathcal{T}=\tm\nh$, while $\,n,\hs r\,$ are the integers in (a) above, 
$\,\vt\,$ is the sub\bd\ of $\,\tm\,$ whose restriction to 
$\,M_\cy\subset M\nh$, for each $\,\cy\in\bs$, is the vertical distribution on 
the $\,\tacb$\pb\ $\,M_\cy$ over $\,Q_\cy$, and $\,\vd\nh=\hs\kerd\pmb\,$ is 
the vertical distribution of the \bp\ $\,\pmb:M\to\bs$. We also set 
$\,\kx(u'\nnh,w)=\beta(u,\zeta)\,$ for any $\,x\in M$ and any vectors 
$\,w\in\txm$, $\,u'\in\vdx=\txmy$ with $\,\cy=\pmb(x)\in\bs$, where 
$\,u\,$ is the $\,\tacb$\inv\ \vf\ tangent to $\,M_\cy$ along the 
$\,\tacb$-orbit of $\,x\,$ and having the value $\,u'$ at $\,x$, while 
$\,\zeta\,$ is the image of $\,w\,$ under the differential at $\,x\,$ of the 
\bp\ $\,M\to Q$.
\vskip3pt
\noindent{\bf Step 2. }{\it We select an arbitrary to\-tal-met\-ric extension 
$\,g\,$ of\/ $\,(\vt,\vd\nnh,\kx)$ restricted to $\,\,U\nh$, where $\,\,U\,$ 
is any fixed nonempty open subset of\/ $\,M$.}
\vskip3pt
\noindent The construction just described gives a null \dn\ $\,\vt$ of 
dimension $\,r\,$ on the $\,n$\diml\ \prd\ $\,(U\nh,g)$. This is clear from 
Lemma~\ref{parfm}, which also implies that such metrics $\,g\,$ are just 
arbitrary sections of some \afb\ over $\,M$.

The reader is again referred to the Appendix for a description of what the 
above steps correspond to in Walker's formulation.

We can now state a co\-or\-di\-nate-free version of Walker's theorem:
\begin{thm}\label{walkr}If\/ $\,g\,$ and\/ $\,\vt\hs$ are obtained as above 
from any prescribed data\/ {\rm(a)} -- {\rm(d)}, then $\,g\,$ is a \prc\ on 
the $\,n$\diml\ \mf\/ $\,\,U\nh$, and\/ $\,\vt\hs$ is a $\,g$-null, $\,g$\prl\ 
distribution of dimension\/ $\,r\hs$ on $\,\,U\nh$.

Conversely, up to an isometry, every \npd\ $\,\vt\,$ on a \prd\/ $\,(M,g)\hs$ 
is, locally, the result of applying the above construction to some data\/ 
{\rm(a)} -- {\rm(d)}. The data themselves are naturally associated with\/ 
$\,g\,$ and\/ $\,\vt$.
\end{thm}
A proof of Theorem~\ref{walkr} is given in the next two sections.

\section{Proof of the first part of Theorem~\ref{walkr}}\label{prfi}
By Lemma~\ref{parfm}, $\,\vt\,$ is $\,g$-null and $\,\vd$ is its 
$\,g$-or\-thog\-o\-nal complement. That 
$\,\vt\,$ is $\,g$\prl\ will be clear if we establish the relation 
$\,\lgl\nabla_{\!w}v,u\rgl=0$, where $\,\nabla\,$ is the \lcc\ of $\,g\,$ and 
$\,\lr\,$ stands for $\,g(\hskip3pt,\hskip1.8pt)$, while $\,v,u,w\,$ are any 
\vf s tangent to $\,M\,$ such that $\,v\,$ is a section of $\,\vt\,$ and 
$\,u\,$ is a section of $\,\vd$. We may further require $\,w\,$ to be 
projectable under both \bd\ projections $\,M\to Q$ and $\,\pmb:M\to\bs$. 
Finally, we may also assume that $\,v\,$ restricted to each $\,\tacb$\pb\ 
space $\,M_\cy$ is an infinitesimal generator of the action of $\,\tacb$, 
while $\,u\,$ restricted to each $\,M_\cy$ is $\,\tacb$\inv. (Locally, such 
$\,w,v,u$ span the \vb s $\,\tm\nnh,\hs\vt\,$ and $\,\vd\nnh$.)

First, $\,[w,v]\,$ is a section of $\,\vt\,$ and $\,[u,w]\,$ is a section of 
$\,\vd$ (from Remark~\ref{liebr} applied to both \bd\ projections), while 
$\,[v,u]=0\,$ by $\,\tacb$-in\-var\-i\-ance of $\,u$. The last three terms in 
(\ref{lcc}) thus all equal zero. 

Our claim will follow if we show that the first three terms in (\ref{lcc}) 
vanish as well. To this end, note that $\,d_w\lgl v,u\rgl=0\,$ since 
$\,\lgl v,u\rgl=0$. Next, $\,d_v\lgl w,u\rgl=0$. Namely, 
$\,\lgl w,u\rgl=\kx(u,w)=\beta(u,\zeta)$, for $\,\kx,\beta,\zeta$ described 
in \hbox{Sec.\ \ref{walk}}, is constant in the direction of $\,v\,$ (and, in 
fact, constant along each leaf of $\,\vt$): at a point 
$\,x\in M_\cy\subset M\,$ we obtain $\,\zeta\,$ as the projection image of 
$\,w(x)$, while $\,u\,$ is $\,\tacb$\inv, so that, due to projectability of 
$\,w$, both $\,u\,$ and $\,\zeta\,$ depend only on the image of $\,x\,$ under 
the \bp\ $\,M\to Q$, rather than $\,x\,$ itself. Finally, 
$\,d_u\lgl w,v\rgl=0\,$ as $\,\lgl w,v\rgl=\xi(\tilde w)\,$ is a function 
$\,\bs\to\bbR\hs$, that is, a function $\,M\to\bbR\,$ constant along $\,\w$. 
Here $\,\xi\,$ is the section of $\,\tab$ corresponding to $\,v\,$ under the 
inclusion $\,\pqb^*(\tab)\subset\mathcal{F}\,$ of \hbox{Sec.\ \ref{walk}}, 
while $\,\tilde w\,$ is the \vf\ on $\,\bs\,$ onto which $\,w\,$ projects; 
therefore, $\,\lgl w,v\rgl=\xi(\tilde w)$, since in \hbox{Sec.\ \ref{walk}} we 
set $\,\hj(\xi,\zeta)=\xi(d\pqb_\pz\zeta)$.

\section{Proof of the second part of Theorem~\ref{walkr}}\label{prse}
For any \npd\ $\,\vt\,$ of dimension $\,\er\hs$ on an $\,n$\diml\ \prd\ 
$\,(M,g)$, the $\,g$-or\-thog\-o\-nal complement $\,\w$ is a parallel \dn\ of 
dimension $\,n-\er$. If the sign pattern of $\,g\,$ has $\,i_-$ minuses and 
$\,i_+$ pluses, it follows that
\begin{equation}
\text{\rm a)}\hskip9pt\er\,\le\,\,\text{\rm min}\hs(i_-,i_+)\hs,
\hskip24pt\text{\rm b)}\hskip9pt\vt\,\subset\,\w\nh,
\hskip24pt\text{\rm c)}\hskip9pt\er\,\le\,n/2\hs.\label{inq}
\end{equation}
In fact, $\,\vt\,$ is null, which gives (\ref{inq}b) and 
$\,\er\le n-\er$, that is, \hbox{(\ref{inq}c),} while (\ref{inq}a) 
follows since, in a pseu\-do\hs-Euclid\-e\-an space with the sign pattern as 
above, $\,i_-$ (or, $\,i_+$) is the maximum dimension of a subspace on which 
the inner product is negative (or, positive) semidefinite.

Every \npd\ $\,\vt\,$ satisfies the curvature relations
\begin{equation}
\text{\rm a)}\hskip5.5ptR(\vt,\nnh\w\hskip-2.8pt,\hh\text{\rm-}\hh,
\hh\text{\rm-}\hh)\nh=\nh0,
\hskip5pt\text{\rm b)}\hskip5.5ptR(\vt,\nnh\vt,\hh\text{\rm-}\hh,
\hh\text{\rm-}\hh)\nh=\nh0,
\hskip5pt\text{\rm c)}
\hskip5.5ptR(\w\hskip-2.8pt,\w\hskip-2.8pt,\vt,\hh\text{\rm-}\hh)\nh
=\nh0,\hskip-13pt
\label{rzo}
\end{equation}
(\ref{rzo}a) meaning that $\,R(v,u,w,w\hh'\hh)=0\,$ whenever $\,v,u,w,w\hh'$ 
are \vf s, $\,v\,$ is a section of $\,\vt$, and $\,u\,$ is a section of 
$\,\w\nnh$. (Similarly for \hbox{(\ref{rzo}b)} and (\ref{rzo}c).) In fact, for 
such $\,v,u,w,w\hh'\nnh$, (\ref{cur}) implies that $\,R(w,w\hh'\hh)v\,$ is a 
section of $\,\vt$, and so it is orthogonal to $\,u$. This proves 
(\ref{rzo}a); (\ref{rzo}a) and (\ref{inq}b) yield (\ref{rzo}b), while 
(\ref{rzo}a) and the first Bianchi identity give (\ref{rzo}c).

We now show how a \npd\ $\,\vt\,$ on a \prd\ $\,(M,g)\,$ gives rise to objects 
(a) -- (d) in \hbox{Sec.\ \ref{walk}}.

First, $\,n\,$ and $\,\er\hs$ are the dimensions of $\,M\,$ and $\,\vt$. By 
(\ref{inq}c), $\,\er\le n/2$.

Being parallel, the distribution $\,\w$ is integrable. Since our discussion is 
local, we will assume, from now on, that $\,M\,$ is the \ts\ of a \bd\ over 
some $\,\er$\diml\ \bmf\ $\,\bs$, whose fibres $\,M_\cy$, $\,\cy\in\bs$, are 
all contractible and coincide with the leaves of $\,\w\nnh$. As $\,\vt$ is 
parallel, the \lcc\ $\,\nabla\,$ induces a connection in the \vb\ obtained by 
restricting $\,\vt\,$ to any given sub\mf\ $\,N$ of $\,M\nh$. In the case 
where $\,N=M_\cy$ is a leaf of $\,\w\nnh$, we have, for each $\,\cy\in\bs$, 
the following conclusion.
\begin{equation}
\begin{array}{ll}
\text{\rm$\tacb\,$ is naturally isomorphic to the space $\,\vy$ of those 
sections}\\
\text{\rm of the restriction of $\,\hs\vt\,$ to $\,M_\cy$ which are parallel 
(along $\,M_\cy$).}
\end{array}\label{nat}
\end{equation}
Instead of establishing (\ref{nat}) directly, we will show that {\it sections 
of\/ $\,\tab$ can be naturally identified with sections of\/ $\,\vt\hs$ 
parallel along $\,\w\nnh$,} using an identification which is clearly 
valuewise, i.e., consists of operators $\,\vy\to\tacb$, $\,\cy\in\bs$. To this 
end, we denote by $\,\pmb\,$ be the \bp\ $\,M\to\bs$. Every \vf\ on $\,\bs\,$ 
is the $\,\pmb\hs$-im\-age $\,(d\pmb)w$ of some $\,\pmb\hs$\pj\ \vf\ $\,w\,$ 
on $\,M$. Let $\,v\,$ now be a section of the \vb\ $\,\vt\,$ over $\,M$, 
parallel in the direction of $\,\w\nnh$. Our identification associates with 
$\,v\,$ the cotangent \vf\ $\,\xi\,$ on $\,\bs$ that sends each \vf\ 
$\,(d\pmb)w\,$ to $\,g(v,w)\,$ treated as a function $\,\bs\to\bbR\hs$. Note 
that $\,\xi\,$ is well defined: two $\,\pmb\hs$\pj\ \vf s $\,w\,$ on $\,M\,$ 
with the same $\,\pmb\hs$-im\-age $\,(d\pmb)w\,$ differ by a section of 
$\,\w\nh=\hs\kerd\pmb$, necessarily orthogonal to $\,v$, so that $\,g(v,w)\,$ 
is the same for both choices of $\,w$. Also, $\,g(v,w):M\to\bbR\,$ actually 
descends to a function $\,\bs\to\bbR\hs$, i.e., is constant along the fibres 
$\,M_\cy$ (leaves of $\,\w$). In fact, $\,d_u[g(v,w)]=0\,$ for any section 
$\,u\,$ of $\,\w$, as $\,\nabla_{\!u}v=0\,$ in view of the assumption about 
$\,v$, and $\,\nabla_{\!u}w=[u,w]+\nabla_{\!w}u$, while $\,[u,w]\,$ (or 
$\,\nabla_{\!w}u$) is a section of $\,\w$ by Remark~\ref{liebr} (or, since 
$\,\w$ is parallel).

Injectivity of the above assignment $\,v\mapsto\xi\,$ is obvious, since 
$\,\pmb\hs$\pj\ \vf s $\,w\,$ span $\,\tm$. Surjectivity of the resulting 
operators $\,\vy\to\tacb\,$ now follows: both spaces have the same dimension, 
as the connections induced by $\,\nabla\,$ in the restrictions of $\,\vt\,$ to 
the leaves $\,M_\cy$ are flat in view of (\ref{rzo}c) (cf.\ (\ref{cur})). This 
proves (\ref{nat}).

Flatness of the induced connections also implies that the leaves of $\,\vt$ 
contained in any given leaf $\,M_\cy$ of $\,\w$ are the fibres of a 
$\,\vy$\pb\ with the \ts\ $\,M_\cy$ over some \bmf\ $\,Q_\cy$. (Here $\,M\,$ 
should be replaced with an open subset, if necessary.) Since each $\,\tacb\,$ 
is identified with $\,\vy$ by (\ref{nat}), we thus obtain the data (c) of 
\hbox{Sec.\ \ref{walk}}.

Next, we define the metric $\,\gm_\cy$ on each $\,Q_\cy$, required by (d) in 
\hbox{Sec.\ \ref{walk}}, so that it assigns the function $\,g(u,u')\,$ to two 
\vf s on $\,Q_\cy$ which are images, under the $\,\tacb$\pb\ projection 
$\,M_\cy\to Q_\cy$, of $\,\tacb$\inv\ \vf s $\,u,u'$ on $\,M_\cy$. Constancy 
of $\,g(u,u')\,$ along the $\,\tacb$-orbits, meaning that $\,d_v[g(u,u')]=0\,$ 
for any section $\,v\,$ of $\,\vt$ defined on $\,M_\cy$ and parallel along 
$\,\w\nnh$, now follows: as $\,v\,$ is $\,\w$\prl\ and $\,u\,$ is 
$\,\tacb$\inv, we have $\,\nabla_{\!u}v=[v,u]=0$, cf.\ (\ref{nat}), so that 
$\,\nabla_{\!v}u=0$. For the same reason, $\,\nabla_{\!v}u'\nh=0$.

Finally, a suitable version of the construction in \hbox{Sec.\ \ref{walk}}, 
applied to the data (a) -- (d) defined above, leads to the original $\,g\,$ 
and $\,\vt$, which is a consequence of how the identification (\ref{nat}) and 
the definition of $\,\gm_\cy$ use $\,g$. The choices of the to\-tal-pair\-ing 
and to\-tal-met\-ric extensions, required in \hbox{Sec.\ \ref{walk}}, are 
provided by $\,g\,$ as well. For instance, $\,\beta\,$ in Step 1 is given by 
$\,\beta(u,\zeta)=g(u,w)$, where $\,u\,$ is a section of $\,\w$ commuting with 
every section $\,v\,$ of $\,\vt\hs$ that is parallel along $\,\w\nnh$, and 
$\,\zeta\,$ is a \vf\ on $\,Q$ (the union of all $\,Q_y$), while $\,w\,$ is 
any \vf\ on $\,M\,$ projectable onto $\,\zeta\,$ under the \bp\ $\,M\to Q$. 
That $\,g(u,w)\,$ depends just on $\,u\,$ and $\,\zeta\,$ (but not on $\,w$) 
is clear: two choices of $\,w\,$ differ by a section of $\,\hs\vt$. Also, 
$\,g(u,w)\,$ is constant in the direction of $\,\hs\vt\hs$ (and so it may be 
treated as a function $\,Q\to M$). Namely, $\,d_v[g(u,w)]=0\,$ for any section 
$\,v\,$ of $\,\vt\hs$ parallel along $\,\w$, which follows as 
$\,\nabla_{\!v}u=\nabla_{\!u}v=0$ (note that $\,[u,v]=0$), while 
$\,\nabla_{\!v}w=[v,w]+\nabla_{\!w}v$, and $\,[v,w]\,$ (or $\,\nabla_{\!w}v$) 
is a section of $\,\vt\,$ by Remark~\ref{liebr} (or, respectively, since 
$\,\vt\,$ is parallel). This completes the proof of Theorem~\ref{walkr}.

\section{The mid\hs-di\-men\-sion\-al case}\label{midi}
For an $\,\er$\diml\ \npd\ $\,\vt\,$ on a \prd\ $\,(M,g)\,$ of dimension 
$\,n=2\er$, the 
discussion in \hbox{Sec.\ \ref{walk}} amounts to nothing new: 
implicitly at least, it is already present in Sec.\ 6 of Walker's original 
paper.$^1$ See also Sec.\ 9 in Ref.~3. (A related global result is 
\hbox{Theorem\ 5} in Ref.~5.) In this section we point out how the 
construction may be simplified when $\,n=2\er$.

Let $\,\vt\,$ and $\,(M,g)\,$ be as above, with $\,n=2\er\ge2$. The relations 
$\,i_-\nh+\,i_+\nh=n\,$ and (\ref{inq}a) imply that $\,g\,$ has the {\it 
neutral\/} sign pattern: $\,i_-\nh=\,i_+\nh=r=n/2$. In (c) -- (d) of 
\hbox{Sec.\ \ref{walk}}, each $\,Q_\cy$ is a $\,0$\diml\ (discrete) 
\mf, and $\,\gm_\cy$ is the ``zero metric'' on $\,Q_\cy$. Also, the choice of 
a to\-tal-pair\-ing extension $\,\beta\,$ in Step 1 of \hbox{Sec.\ \ref{walk}} 
is now unique: the \afb\ having $\,\beta\,$ as a section is of fibre dimension 
$\,0$. The construction in \hbox{Sec.\ \ref{walk}} can therefore be rephrased 
as follows. Given
\begin{enumerate}
  \def\theenumi{{\rm\alph{enumi}}}
\item an even integer $\,n\ge2$,
\item a \mf\ $\,\bs\,$ of dimension $\,r=n/2$,
\item an \afb\ over $\,\bs\,$ with some \ts\ $\,M\nh$, for which $\,\tab\,$ 
is the associated \vb\ (\hbox{Sec.\ \ref{afpb}}),
\end{enumerate}
we define a partial metric $\,(\vt,\vd\nnh,\kx)\,$ in the tangent \bd\ 
$\,\tm\,$ by choosing $\,\vt=\vd$ to be the vertical distribution 
$\,\hs\kerd \pmb\,$ for the \bp\ $\,\pmb:M\to\bs$, and setting 
$\,\kx(\xi,w)=\xi(d\pmb_xw)\,$ for any $\,x\in M$, 
$\,\xi\in\vt_x=\tacb$, where $\,\cy=\pmb(x)$, and $\,w\in \txm$. Selecting any 
to\-tal-met\-ric extension $\,g\,$ of $\,(\vt,\vd\nnh,\kx)\,$ on a fixed 
nonempty open set $\,\,U\subset M$, we now obtain an $\,n$\diml\ 
\prd\ $\,(U\nh,g)\,$ on which $\,\vt\,$ is a $\,g$-null, $\,g$\prl\ 
distribution of dimension $\,\er=n/2$.

Conversely, up to an isometry, every \npd\ $\,\vt$ of dimension 
$\,\er\ge1\,$ on a \prd\ $\,(M,g)\,$ with $\,\dim M=2\er\hs$ arises, locally, 
from the above construction applied to some data (a) -- (c), themselves 
naturally determined by $\,g\,$ and $\,\vt$.

\section*{Acknowledgments}
The authors wish to thank Zbigniew Olszak for helpful comments.

\section*{Appendix: Walker's original statement}
Walker stated his classification result as follows.$^1$
\vskip2pt
T{\eightrm HEOREM} 1. {\it A canonical form for the general\/ $\,V_n$ of class 
$\,C^\infty$ $(\hskip-.9pt\text{\it or\ }\,C^{\hs\omega})$ admitting a 
parallel null\/ $\,r$-plane is given by the fundamental tensor
\vskip2pt
\centerline{$
(g_{ij})\,=\,\left(\begin{matrix}O&O&I\\O&A&H\\I&\hskip1ptH'&\hskip-1ptB
\end{matrix}\right)$}
\vskip5pt
\noindent where $\,I\hs$ is the unit\/ $\,\hs r\times\hs r\,$ matrix and\/ 
$\,A,\,B,\,H,\,H'$ are matrix functions of the coordinates, of the same class 
as $\,V_n$, satisfying the following conditions but otherwise arbitrary\/{\rm:}
\begin{enumerate}
  \def\theenumi{{\rm\roman{enumi}}}
\item $A\,$ and\/ $\,B\,$ are symmetric, $\,A\,$ is of order 
$\,(n-2r)\times(n-2r)\,$ and nonsingular, $\,B\,$ is of order $\,r\times r$, 
$\,H\,$ is of order $\,(n-2r)\times r$, and\/ $\,H'$ is the transpose of\/ 
$\,H$.
\item $A\,$ and\hskip3pt$\,H\hskip3pt(\text{\it and\ therefore\ }H'\hs)$ are 
independent of the coordinates $\,x^1,\dots,x^r\hskip-2.5pt$.
\end{enumerate}

A basis for the parallel null\/ $\,r$-plane is the set of vectors 
$\,\delta_1^i,\delta_2^i,\dots,\delta_r^i$.}

\vskip7pt
\noindent Here is how the coordinates and matrix functions appearing above 
correspond to the objects used for the construction in \hbox{Sec.\ 
\ref{walk}}. Walker's coordinates $\,x^i\nh$, $\,i=1,\dots,n$, serve as a 
coordinate system for the \mf\ $\,M\,$ of \hbox{Sec.\ \ref{walk}}. 
Coordinates for other manifolds appearing in \hbox{Sec.\ \ref{walk}} are 
obtained from $\,x^i$ by restricting the range of the index $\,i$, to 
$\,i>n-\er\hs$ (for $\,\bs$), $\,i>\er\hs$ (for $\,Q$), $\,i\le n-\er\hs$ (for 
each $\,M_\cy$) and $\,\er<i\le n-\er\hs$ (for each $\,Q_\cy$). The center 
submatrix $\,A\,$ in Walker's matrix corresponds to the family $\,\gm_\cy$, 
$\,\cy\in\bs$, of \prc s ((d) in \hbox{Sec.\ \ref{walk}}) and, consequently, 
also to the formula for $\,\hj(u,\psi)$, while the last two matrices 
$\,O\hskip4ptI\,$ in the first row represent the definition of 
$\,\hj(\xi,\zeta)$. The Walker-matrix counterpart of the extension $\,\beta\,$ 
chosen in Step 1 is the $\,(n-\er)\times(n-\er)\,$ submatrix with the rows 
$\,O\hskip4ptI\,$ and $\,A\hskip4ptH$, so that the freedom in choosing 
$\,\beta\,$ amounts to arbitrariness in the selection of $\,H\,$ (and $\,H\,$ 
is independent of the coordinates $\,x^i\nh$, $\,i=1,\dots,\er$, which 
translates into the fact that $\,\beta\,$ is a morphism of \vb s over the \mf\ 
$\,Q\,$ with the coordinates $\,x^i\nh$, $\,i>\er$). Once chosen, $\,\beta\,$ 
is used in \hbox{Sec.\ \ref{walk}} to define $\,\vt,\vd$ and $\,\kx$. In 
terms of Walker's coordinates and matrix functions, $\,\vt\,$ (or, $\,\vd$) 
is spanned by the $\,x^i$ coordinate directions with $\,i\le\er\hs$ (or, 
respectively, $\,i\le n-\er$), while the analog of $\,\kx\,$ is the 
$\,(n-\er)\times n\,$ submatrix with the rows $\,O\hskip4ptO\hskip4ptI\,$ and 
$\,O\hskip3.4ptA\hskip4ptH$. Finally, the extension in Step 2 is nothing else 
than augmenting this last submatrix by a third row, 
$\,I\hskip4ptH'\hskip1.4ptB$, in which $\,B\,$ is completely arbitrary.

\newpage
\hskip-43pt\parbox[l]{391pt}{
{\eightpoint
\begin{enumerate}
\item[\ $^1$\hskip-5pt]Walker,\hskip2.3ptA.\hskip1.6ptG., ``Canonical form 
for a Riemannian space with a parallel field of null planes,'' Quart.\ J.\ 
Math.\ Oxford (2) {\eb1} (1950), 69--79.
\item[\ $^2$\hskip-5pt]Walker,\hskip2.3ptA.\hskip1.6ptG., ``Canonical forms. 
II. Parallel partially null planes,'' Quart.\ J.\ Math.\ Oxford (2) {\eb1} 
(1950), 147--152; ``Riemann extensions of \hbox{non\hs-}Riem\-ann\-i\-an 
spaces,'' {\sii Con\-ve\-gno In\-ter\-na\-zio\-na\-le di Ge\-o\-me\-tria 
Differenziale, Italia, 1953} (E\-di\-zio\-ni Cre\-mo\-ne\-se, Ro\-ma, 1954, 
pp.\ 64--70).
\item[\ $^3$\hskip-5pt]Patterson, \hskip2.3ptE.\hskip1.6ptM., and 
Walker,\hskip2.3ptA.\hskip1.6ptG., ``Riemann extensions,'' Quart.\ J.\ Math.\ 
Oxford (2) {\eb3} (1952), 19--28.
\item[\ $^4$\hskip-5pt]Patterson, \hskip2.3ptE.\hskip1.6ptM., ``Riemann 
extensions which have K\"ahler metrics,'' Proc.\ Roy.\ Soc.\ Edinburgh, Sect.\ 
A. {\eb64} (1954), 113--126.
\item[\ $^5$\hskip-5pt]Robertson,\hskip2.3ptS.\hskip1.6ptA. and 
Furness,\hskip2.3ptP.\hskip1.6ptM.\hskip1.6ptD., ``Parallel framings and 
foliations on pseu\-\hbox{do\hs-}Riem\-ann\-i\-an manifolds,'' J.\ Diff.\ 
Geometry {\eb9} (1974), 409-422.
\item[\ $^6$\hskip-5pt]Roter,\hskip2.3ptW., ``On conformally symmetric 
Ric\-ci-re\-cur\-rent spaces,'' Colloq.\ Math. {\eb31} (1974), 87--96.
\item[\ $^7$\hskip-5pt]Farran, \hskip2.3ptH., ``Foliated 
pseu\-\hbox{do\hs-}Riem\-ann\-i\-an and symplectic manifolds,'' Progr.\ Math.\ 
(Al\-lah\-abad) {\eb13} (1979), 59--64.
\item[$^8$\hskip-5pt]Egorov,\hskip2.3ptA.\hskip1.6ptI., ``Parallel 
isotropic bivector area in $\,V\sb{4}$'' (in Russian), Gravitatsiya i Teor.\ 
Ot\-no\-si\-\hbox{tel\hskip-.2pt\vbox{\hbox{$^{'}$}\vskip-2.5pt}}\hskip-1.8ptno\-sti {\eb20} (1983), 20--24.
\item[$^9$\hskip-5pt]Thompson,\hskip2.3ptG., ``Normal form of a metric 
admitting a parallel field of planes,'' J.\ Math.\ Phys. {\eb33} (1992), 
4008--4010.
\item[$^{10}$\hskip-5pt]Ghanam,\hskip2.3ptR. and Thompson,\hskip2.3ptG., 
``The holonomy Lie algebras of neutral metrics in dimension four,'' J.\ Math.\ 
Phys. {\eb42} (2001), 2266--2284.
\item[$^{11}$\hskip-5pt]Matsushita,\hskip2.3ptY., ``Four-di\-men\-sion\-al 
Walker metrics and symplectic structures,'' J.\ Geom.\ Phys. {\eb52} (2004), 
89--99; ``Walker $\,4$-man\-i\-folds with proper almost complex structures,'' 
{\sii ibid.}, {\eb55} (2005), 385--398.
\item[$^{12}$\hskip-5pt]Chaichi,\hskip2.3ptM., 
Garc\'\i a\hh-R\'\i o,\hskip2.3ptE., and Matsushita,\hskip2.3ptY., 
``Curvature properties of four-di\-men\-sion\-al Walker metrics,'' Classical 
Quantum Gravity {\eb22} (2005), 559--577.
\item[$^{13}$\hskip-5pt]Chaichi,\hskip2.3ptM., 
\hskip-1ptGarc\'\i a\hh-R\'\i o,\hskip2.3ptE., and 
V\hskip-.9pt\'azquez-Abal,\hskip2.3ptM.\hskip1.6ptE., 
``Three-di\-men\-sion\-al Lo\-rentz manifolds admitting a parallel null vector 
field,'' J.\ Phys. {\eb A 38} (2005), 841--850.
\item[$^{14}$\hskip-5pt]Matsushita,\hskip2.3ptY., Haze,\hskip2.3ptS., and 
Law,\hskip2.3ptP.\hskip1.6ptR., ``Almost K\"ah\-ler-Ein\-stein structures 
on $\,8$-di\-men\-sion\-al Walker manifolds,'' Monatsh.\ Math., in press.
\item[$^{15}$\hskip-5pt]Ovando,\hskip2.3ptG., ``Invariant pseudo K\"ahler 
metrics in dimension four,'' math.DG/0410232.
\item[$^{16}$\hskip-5pt]Derdzinski,\hskip2.3ptA. and Roter,\hskip2.3ptW., 
``Projectively flat surfaces, null parallel distributions, and conformally 
symmetric manifolds,'' preprint (URL 
http:/\hskip-1.5pt/www.math.ohio\hh-state.edu/\~{}\hskip-.8ptandrzej).
\item[$^{17}$\hskip-5pt]Kobayashi,\hskip2.3ptS. and Nomizu,\hskip2.3ptK., 
{\sii Foundations of \hskip.8ptDifferential Geometry} (Wiley, New York, 1963), 
\hbox{Vol.\ I,} p.\ 160.
\end{enumerate}
}
}
\end{document}